\documentclass[12pt,british]{article}

\usepackage{babel}
\usepackage{amsmath,amsfonts,amssymb,amsthm}
\usepackage[latin1]{inputenc}
\usepackage[T1]{fontenc}
\usepackage{enumerate}

\makeatletter

\renewcommand{\p@enumii}{}

\def\@enum@{\list{\csname label\@enumctr\endcsname}%
           {\usecounter{\@enumctr}\def\makelabel##1{
\normalfont\ignorespaces\emph{{##1}~}}
\setlength{\labelsep}{3pt}
\setlength{\parsep}{0pt}
\setlength{\itemsep}{0pt}
\setlength{\leftmargin}{0pt}
\setlength{\labelwidth}{0pt}
\setlength{\listparindent}{\parindent}
\setlength{\itemsep}{0pt}
\setlength{\itemindent}{0pt}
\topsep=3pt plus 1pt minus 1 pt}}

\makeatother

\renewcommand{\epsilon}{\ensuremath{\varepsilon}}
\renewcommand{\to}{\ensuremath{\rightarrow}}

\newcommand{\R}{\ensuremath{\mathbb R}}
\newcommand{\N}{\ensuremath{\mathbb N}}
\newcommand{\Z}{\ensuremath{\mathbb Z}}
\newcommand{\dt}{\ensuremath{\mathbb D}^{2}}
\newcommand{\St}[1][2]{\ensuremath{\mathbb S}^{#1}}

\DeclareRobustCommand*{\up}[1]{\textsuperscript{#1}}

\newcommand{\ft}[1][n]{\ensuremath{T_{#1}}}
\renewcommand{\ker}[1]{\ensuremath{\operatorname{\text{Ker}}\left({#1}\right)}}

\makeatletter
\def\@map#1#2[#3]{\mbox{$#1 \colon\thinspace #2 \to #3$}}
\def\map#1#2{\@ifnextchar [{\@map{#1}{#2}}{\@map{#1}{#2}[#2]}}
\makeatother

\newcommand{\brak}[1]{\ensuremath{\left\{ #1 \right\}}}
\newcommand{\ang}[1]{\ensuremath{\langle #1\rangle}}

\newcommand{\si}[2][{}]{\ensuremath{\sigma_{#2}^{#1}}}

\newtheoremstyle{theoremm}{}{}{\itshape}{}{\scshape}{.}{ }{}

\theoremstyle{theoremm}
\newtheorem*{thm}{Theorem}

\newtheoremstyle{remarkk}{}{}{}{}{\scshape}{.}{ }{}

\theoremstyle{remarkk}
\newtheorem*{rem}{Remark}

\begin{document}

\title{The quaternion group as a subgroup of the sphere braid groups}

\author{DACIBERG~LIMA~GON\c{C}ALVES\\
Departamento de Matem\'atica - IME-USP,\\
Caixa Postal~\textup{66281}~-~Ag.~Cidade de S\~ao Paulo,\\ CEP:~\textup{05311-970} - S\~ao Paulo - SP - Brazil.\\ e-mail\textup{:~\texttt{dlgoncal@ime.usp.br}}\vspace*{4mm}\\
JOHN~GUASCHI\\
Laboratoire de Math\'ematiques Emile Picard, UMR CNRS~\textup{5580},\\ UFR-MIG, Universit\'e Toulouse~III, \textup{31062}~Toulouse Cedex~\textup{9}, France.\\
e-mail\textup{:~\texttt{guaschi@picard.ups-tlse.fr}}}

\date{\today}

\begingroup
\renewcommand{\thefootnote}{}
\endgroup 


\maketitle

\begin{abstract}\noindent
\emph{Let $n\geq 3$. We prove that the quaternion group of order $8$ is realised as a subgroup of the sphere braid group $B_n(\St)$ if and only if $n$ is even. If $n$ is divisible by $4$ then the commutator subgroup of $B_n(\St)$ contains such a subgroup. Further, for all $n\geq 3$, $B_n(\St)$ contains a subgroup isomorphic to the dicyclic group of order $4n$.}
\end{abstract}

The braid groups $B_n$ of the plane were introduced by E.~Artin
in~1925~\cite{A1,A2}, and were generalised by Fox to braid groups of
arbitrary topological spaces using the notion of configuration
space~\cite{FoN}. Van Buskirk showed that the braid groups of a
compact connected surface $M$ possess torsion elements if and only if
$M$ is the sphere $\St$ or the real projective plane $\R
P^2$~\cite{vB}. Let us recall briefly some of the properties of the
braid groups of the sphere~\cite{FvB,GVB,vB}.

If $\dt\subseteq \St$ is a topological disc, there is a group
homomorphism $\map {\iota}{B_n}[B_n(\St)]$ induced by the
inclusion. If $\beta\in B_n$ then we shall denote its image
$\iota(\beta)$ simply by $\beta$. Then $B_n(\St)$ is
generated by $\sigma_1,\ldots,\sigma_{n-1}$ which are subject to the
following relations:
\begin{equation*}\label{eq:presnbns}
\begin{gathered}
\text{$\si{i}\si{j}=\si{j}\si{i}$ if $\lvert i-j\rvert\geq 2$ and $1\leq i,j\leq n-1$}\\
\text{$\si{i}\si{i+1}\si{i}=\si{i+1}\si{i}\si{i+1}$ for all $1\leq i\leq n-2$, and}\quad\\
\text{$\si 1\cdots \si {n-2}\si [2]{n-1}\si {n-2}\cdots\si 1=1$.}
\end{gathered}
\end{equation*}
Consequently, $B_n(\St)$ is a quotient of $B_n$. The first three sphere braid groups are finite: $B_1(\St)$ is trivial, $B_2(\St)$ is cyclic of order~$2$, and $B_3(\St)$ is a $\text{ZS}$-metacyclic group (a group whose Sylow subgroups, commutator subgroup and commutator quotient group are all cyclic) of order~$12$. The Abelianisation of $B_n(\St)$ is isomorphic to the cyclic group $\Z_{2(n-1)}$. The kernel of the associated projection $\map{\xi}{B_n(\St)}[\Z_{2(n-1)}]$ (which is defined by $\xi(\sigma_i)=\overline{1}$ for all $1\leq i\leq n-1$) is the commutator subgroup $\Gamma_2\left(B_n(\St) \right)$. If $w\in B_n(\St)$ then $\xi(w)$ is  the exponent sum (relative to the $\sigma_i$) of $w$ modulo $2(n-1)$. 

The torsion elements of the braid groups of $\St$ and $\R P^2$ were
classified by Murasugi~\cite{M}: if $M=\St$ and $n\geq 3$, they are
all conjugates of powers of the three elements $\alpha_0=
\sigma_1\cdots \sigma_{n-2} \sigma_{n-1}$ (which is of order $2n$),
$\alpha_1=\sigma_1\cdots \sigma_{n-2} \sigma_{n-1}^2$ (of order
$2(n-1)$) and $\alpha_2=\sigma_1\cdots \sigma_{n-3} \sigma_{n-2}^2$
(of order $2(n-2)$) which are respectively $n\up{th}$, $(n-1)\up{th}$
and $(n-2)\up{th}$ roots of $\ft$, where $\ft$ is the so-called `full
twist' of $B_n(\St)$, defined by $\ft=
(\sigma_1\cdots\sigma_{n-1})^n$. If $n\geq 3$, $\ft[n]$ is the unique
element of $B_n(\St)$ of order $2$ and generates the centre of
$B_n(\St)$. In~\cite{GG4}, we showed that $B_n(\St)$ is generated by
$\alpha_0$ and $\alpha_1$.

For $n\geq 4$, $B_n(\St)$ is infinite. It is an interesting question
as to which finite groups are realised as subgroups of $B_n(\St)$
(apart of course from the cyclic groups $\ang{\alpha_i}$).
In~\cite{GG4}, we proved that $B_n(\St)$ contains an isomorphic copy
of the finite group $B_3(\St)$ of order $12$ if and only if
$n\not\equiv 1 \bmod 3$.  The quaternion group $\mathcal{Q}_8$ of
order $8$ appears in the study of braid groups of non-orientable
surfaces, being isomorphic to the $2$-string pure braid group $P_2(\R
P^2)$. Further, since the projection $F_3(\R P^2)\to F_2(\R P^2)$ of
configuration spaces of $\R P^2$ onto the first two coordinates admits
a section~\cite{vB}, it follows using the Fadell-Neuwirth short exact
sequence that $P_3(\R P^2)$ is a semi-direct product of a free group
of rank~$2$ by $\mathcal{Q}_8$~\cite{GG3}.

While studying the lower central and derived series of the sphere braid groups, we showed that $\Gamma_2\left(B_4(\St) \right)$ is isomorphic to a semi-direct product of $\mathcal{Q}_8$ by a free group of rank $2$~\cite{GG5}. After having proved this result, we noticed that the question of the realisation of $\mathcal{Q}_8$ as a subgroup of $B_n(\St)$ was explicitly posed by R.~Brown in connection with the fact that the fundamental group of $\operatorname{SO}(3)$ is isomorphic to $\Z_2$~\cite{DD}. In this paper, we give a complete answer to this question:
\begin{thm}
Let $n\in\N$, $n\geq 3$.
\begin{enumerate}
\item\label{it:q8div4} If $n$ is a multiple of $4$ then $\Gamma_2\left(B_n(\St) \right)$ contains a subgroup isomorphic to $\mathcal{Q}_8$.
\item\label{it:q8even} If $n$ is an odd multiple of $2$ then $B_n(\St)$ contains a subgroup isomorphic to $\mathcal{Q}_8$.
\item\label{it:q8odd} If $n$ is odd then $B_n(\St)$ contains no subgroup isomorphic to $\mathcal{Q}_8$.
\end{enumerate}
\end{thm}

\begin{proof}[\textsc{Proof}]
We first suppose that $n$ is even, so that $n=2m$ with $m\in\N$. Let $H$ be the subgroup of $B_{2m}(\St)$ generated by $x$ and $y$, where
\begin{align*}
x=& (\sigma_1\cdots \sigma_{2m-1}) (\sigma_1\cdots \sigma_{2m-2}) \cdots (\sigma_1 \sigma_2) \sigma_1,\\
y=& (\sigma_1\cdots \sigma_{m-1}) (\sigma_1\cdots \sigma_{m-2}) \cdots (\sigma_1 \sigma_2) \sigma_1\cdot \sigma_{2m-1}^{-1} (\sigma_{2m-2}^{-1} \sigma_{2m-1}^{-1})\cdots\\
& \cdots (\sigma_{m+2}^{-1}\cdots \sigma_{2m-1}^{-1}) (\sigma_{m+1}^{-1}\cdots \sigma_{2m-1}^{-1}).
\end{align*}
Geometrically, $x$ may be interpreted as the `half twist' or Garside element of $B_{2m}$~\cite{Bi}. Further, $y$ may be considered as the commuting product of the positive half twist of the first $m$ strings with the negative half twist of the last $m$ strings. Then $x^2=\ft[2m]$ and $y^2= (\sigma_1 \ldots \sigma_{m-1})^m (\sigma_{m+1}^{-1}\ldots \sigma_{2m-1}^{-1})^m= \ft[2m]$ in $B_{2m}(\St)$ (cf.~\cite{FvB,GVB}). It is well known that $x\sigma_i x^{-1}=\sigma_{2m-i}$ in $B_{2m}$~\cite{Bi}, and thus in $B_{2m}(\St)$, from which we obtain $xyx^{-1}=y^{-1}$. Hence $H$ is isomorphic to a quotient of $\mathcal{Q}_8$. But $x$ is of order $4$, and the induced permutation of $y$ on the symmetric group $S_{2m}$ is different from that of the elements of $\ang{x}$. It follows that $H$ contains the five distinct elements of $\ang{x}\cup \brak{y}$, and so $H\cong \mathcal{Q}_8$. If $m$ is even then $x,y\in \ker{\xi}$, and thus $H\subseteq \Gamma_2\left(B_n(\St) \right)$. This proves parts~(\ref{it:q8div4}) and~(\ref{it:q8even}).

To prove part~(\ref{it:q8odd}), suppose that $n$ is odd, and suppose that $x,y\in B_n(\St)$ generate a subgroup $H$ isomorphic to $\mathcal{Q}_8$, so $x^2=y^2$ and $xyx^{-1}=y^{-1}$. In particular, $x$ and $y$ are of order $4$, and thus are conjugates of $\alpha_1^{\pm (n-1)/2}$ by Murasugi's classification. By considering a conjugate of $H$ if necessary, we may suppose that $x=\alpha_1^{\epsilon_1 (n-1)/2}$ and $y=w\alpha_1^{\epsilon_2 (n-1)/2}w^{-1}$, where $w\in B_n(\St)$ and $\epsilon_1, \epsilon_2\in\brak{1,-1}$. Replacing $x$ by $x^{-1}$ if necessary, we may suppose further that $\epsilon_1=-\epsilon_2$. Thus $xy=\left[\alpha_1^{\epsilon_1 (n-1)/2} ,w\right]$, and is of exponent sum zero. On the other hand, $xy$ is an element of $H$ of order $4$, and so is conjugate to $\alpha_1^{\pm (n-1)/2}$ by Murasugi's classification. But $\xi\left(\alpha_1^{\pm (n-1)/2}\right)=\pm \frac{n(n-1)}{2}$, which is non zero modulo $2(n-1)$. This yields a contradiction, and proves part~(\ref{it:q8odd}).
\end{proof}

\begin{rem}
Let $n\geq 3$. Using techniques similar to those of the proof of the~Theorem, one may show that the subgroup of $B_n(\St)$ generated by $\sigma_1\cdots \sigma_{n-1}$ and the half twist $x$ is isomorphic to the dicyclic group of order $4n$. In particular, if $n$ is a power of two then $B_n(\St)$ contains a subgroup isomorphic to the generalised quaternion group of order $4n$. Further investigation into the finite subgroups of $B_n(\St)$ and $B_n(\R P^2)$ will appear elsewhere.
\end{rem}

\end{document}